\documentclass{ifacconf}

\usepackage{graphicx}      
\usepackage{natbib}        

\usepackage{balance}

\usepackage{amsmath} 

\usepackage{color}

\usepackage{url}

\DeclareFontFamily{U}{msb}{}
\DeclareFontShape{U}{msb}{m}{n}{ <5> <6> <7> <8> <9> gen * msbm
<10> <10.95> <12> <14.4> <17.28> <20.74> <24.88> msbm10}{}
\DeclareSymbolFont{AMSb}{U}{msb}{m}{n}
\DeclareMathSymbol{\inertia}{\mathalpha}{AMSb}{'111}
\DeclareMathSymbol{\Reals}{\mathalpha}{AMSb}{'122}
\DeclareMathSymbol{\Naturals}{\mathalpha}{AMSb}{'116}
\DeclareMathSymbol{\Knumbers}{\mathalpha}{AMSb}{'113}
\DeclareMathSymbol{\Rationals}{\mathalpha}{AMSb}{'121}
\DeclareSymbolFont{AMSb}{U}{msb}{m}{n}
\DeclareMathSymbol{\setB}{\mathalpha}{AMSb}{'102}
\DeclareMathSymbol{\setC}{\mathalpha}{AMSb}{'103}
\DeclareMathSymbol{\setD}{\mathalpha}{AMSb}{'104}
\DeclareMathSymbol{\setE}{\mathalpha}{AMSb}{'105}
\DeclareMathSymbol{\setF}{\mathalpha}{AMSb}{'106}
\DeclareMathSymbol{\setI}{\mathalpha}{AMSb}{'111}
\DeclareMathSymbol{\setK}{\mathalpha}{AMSb}{'113}
\DeclareMathSymbol{\setM}{\mathalpha}{AMSb}{'115}
\DeclareMathSymbol{\setN}{\mathalpha}{AMSb}{'116}
\DeclareMathSymbol{\setP}{\mathalpha}{AMSb}{'120}
\DeclareMathSymbol{\setQ}{\mathalpha}{AMSb}{'121}
\DeclareMathSymbol{\setR}{\mathalpha}{AMSb}{'122}
\DeclareMathSymbol{\setS}{\mathalpha}{AMSb}{'123}

\begin{document}
\begin{frontmatter}

\title{Power Network Regulation Benchmark for Switched-Mode Optimal Control}

\author[First]{Timothy M. Caldwell} 
\author[Second]{Todd D. Murphey} 
\address[First]{University of Colorado, 
   Boulder, CO 80309 USA (e-mail: caldwelt@colorado.edu).}
\address[Second]{Northwestern University, 
   Evanston, IL 60208 USA (e-mail: t-murphey@northwestern.edu)}

\begin{abstract}
Power network regulation is presented as a benchmark problem for assessing and developing switched-mode optimal control approaches like mode scheduling, sliding window scheduling and modal design. Power network evolution modeled by the swing equations and coupled with controllable switching components is a nonlinear, high-dimensional problem. The proposed benchmark problem is the 54 generator IEEE 118 Bus Test Case composed of 106 states. Open questions include scalability in state and number of modes of operation, as well as real-time implementation, reliability, hysteresis, and timing constraints. Can the entire North American power network be regulated? Can every transmission line have independent switching control authority? 
\end{abstract}

\begin{keyword}
switched-mode optimal control, power network regulation, benchmark examples, mode scheduling
\end{keyword}

\end{frontmatter}

\section{Introduction}
\begin{figure*}
\centering
\includegraphics[width = 400pt]{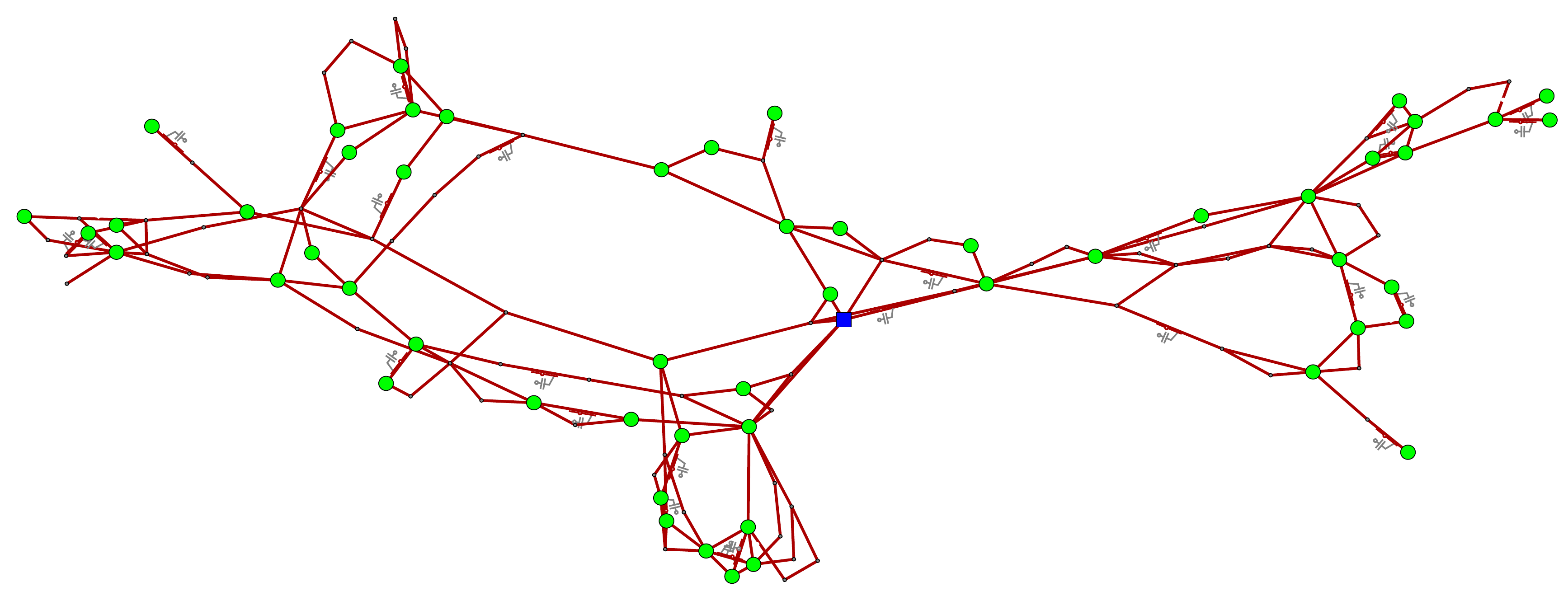}
\caption{Representation of the IEEE 118 Bus Test Case.  The network is composed of 118 buses, 186 lines, and 54 generators.  The generators are depicted by green circles and the reference generator is depicted by the blue square.  The location of the 26 capacitor are shown.}
\label{fig-case118_Graph}
\end{figure*}

This paper presents a benchmark problem for switched-mode systems regulation of a power network. The power network's dynamics are dictated by the swing equations \citep{grainger_stevenson}, which are sufficiently rich to both assess and develop new switched-mode system control strategies. Power networks have a large number of generating nodes that operate as coupled oscillators for which conditions for synchronization have been studied \cite{dorfler_chertkov_bullo}. We investigate the problem of active control to respond to a disturbance. For benchmark, we propose study of the 54 generator IEEE 118 Bus Test Case \citep{christie} depicted in Fig.~\ref{fig-case118_Graph}. The Test Case has 106 states composed of each generator's rotor angle and angular velocity.

The swing equations are nonlinear differential equations, which, when coupled with controllable switching components, are switched-mode.  The proposed benchmark problems are practical for assessing and developing switched-mode system control and design strategies like mode scheduling, sliding window scheduling, and modal design as proposed in this paper. 

Mode scheduling is an integer-constrained optimal control problem. Solving the benchmark problem through direct methods \citep{alamir_attia,lincoln_bernhardsson, lincoln_rantzer} is an optimization with dimension in the range of $10^{10}$ to $10^{11}$. Direct methods require a fine timing discretization in order to capture equivalent fast switching solutions of indirect methods like the projection-based mode scheduling method we presented and analyzed in \citep{caldwell_murphey_nahs}. Implemented in a sliding window framework, the projection-based mode scheduling method successfully computes a solution at roughly twenty times slower than real-time on a standard desktop machine.  

Benchmark questions include both scalability in state (can the entire U.S. network be regulated?), scalability in switching signals (can every transmission line have independent switching control authority?), reliability, hysteresis, timing constraints, et cetera.

Numerical methods vary for solving the optimal mode schedule. Many discretize in time \citep{alamir_attia,lincoln_bernhardsson, lincoln_rantzer} while others in space, \citep{bemporad_borrelli_morari, shaikh_caines_surf, caines_shaikh_2005, giua_seatzu_vandermee_1,giua_seatzu_vandermee_2, hedlund_rantzer, schild_ding_egerstedt_lunze, seatzu_corona_etal}, but often such approaches result in a combinatoric search \citep{shaikh_caines_comb, caines_shaikh_2005, gorges_izak_liu, lincoln_bernhardsson, lincoln_rantzer, seatzu_corona_etal}.  An alternative is embedding methods \citep{bengea_decarlo} which transforms the space of controls to boundary constrained functions and relies on the fact that the set of switched system trajectories is dense in the set of embedded trajectories.  A final alternative taken in \citep{axelsson_wardi_et_al, caldwell_murphey_nahs, egerstedt_wardi_axelsson, gonzalez_vasudevan_etal, kawashima_et_al, wardi_egerstedt_acc12,wardi_egerstedt_hale} pursues variational methods which repeatedly vary the schedule to reduce the cost. The projection-based mode scheduling algorithm \citep{caldwell_murphey_nahs} is a member of this final group.  

This paper is organized as follows: Section~\ref{sec-benchmark} introduces the power network regulation benchmark along with the IEEE 118 Bus Test Case. Section~\ref{sec-problems} presents and provides results for mode scheduling, sliding window scheduling, and modal design of the benchmark problem. Section~\ref{sec-discussion} discusses switched-mode optimal control extensions and future advancements using the benchmark problem.

\section{The Benchmark System}
\label{sec-benchmark}
Due to the complex interconnectedness of multimachine power networks, it is unclear how to actively respond to a disturbance.  The solution we propose is to compute a schedule for physical switches that connect and disconnect capacitors from the network so that system performance improves. A power network is often modeled as a synchronous machine where the dynamics are given by the swing equations \citep{grainger_stevenson}.  The swing equations are second-order nonlinear differential equations which dictate the evolution of each generator's rotor angle.  The rotors are assumed to be spinning at a constant frequency\textemdash e.g. 60 Hz\textemdash but each rotor's relative phase may not be constant.  The evolution of a single rotor is dictated by the difference of its relative phase with its neighboring generators as well as the admittance of the adjacent transmission lines and buses. Through switching capacitors, the transmission lines' admittance value switches, effectively splitting the system dynamics into distinct operating modes dependent on the position of the switches. The only control authority we impose is through the switches.

Let $\delta_i(t)$ be the rotor phase of generator $i$ relative to a reference generator, generator 0.  The evolution of the $i^\textrm{th}$ rotor is dictated by the difference of the mechanical power input with the electrical power output:
\begin{equation}
\frac{2H_i}{\omega_s}\ddot{\delta}_i = P_{m,i} - P_{e,i}
\label{eq-pn_dynamics}
\end{equation}
where $\omega_s$, in $rads/s$ is the synchronous speed and $H_i$ is a normalized inertial constant so that the mechanical power $P_{m,i}$ and the electrical power $P_{e,i}$ are in per unit.  The terms $\omega_s$, $H_i$ and $P_{m,i}$ are assumed to be constant for the short time horizon for which the disturbance and resolution occurs. The electrical power output of generator $i$, $P_{e,i}$, depends on the difference of its rotor's relative phase with the neighboring generators' as well as the admittance of the adjacent lines and buses:
\begin{equation}
P_{e,i} = |E_i|^2G_{ii}+\sum_{j\neq i}|E_i||E_j||Y_{ij}|\cos(\delta_i-\delta_j-\psi_{ij})
\label{eq-powerout}
\end{equation}
where $E_i$ is the transient internal voltage, $G_{ii}$ is the real part of the $ii^\textrm{th}$ component of the bus admittance matrix, $Y_{ij}$ is the $ij^\textrm{th}$ component of the bus admittance matrix, and $\psi_{ij}$ is the angle of the $ij^\textrm{th}$ component of the bus admittance matrix.  

The power network state is composed of the rotor angles and the rotor angular velocities $x(t) = [\delta_1(t),\ldots,\delta_{n/2}(t),$ $\dot{\delta}_1(t),\ldots,\dot{\delta}_{n/2}(t)]^T$, where $n/2$ is the number of generators. Through Eq.~(\ref{eq-pn_dynamics}), the dynamics can be written in the general form
\[
\dot{x}(t) = f(x(t))\textrm{, }x(0) = x_0
\]
with initial condition $x_0$. A disturbance is modeled as an initial perturbation from a known steady state.

\subsection{Switching Control}
To impose a control, we switch the line admittances $Y_{ij}$ between set values by connecting and disconnecting capacitors from transmission lines placed throughout the network. With such a control authority, the power network is an autonomous switched system with distinct operating modes. Depending on the state of the switched capacitor banks at any given time, the system evolves according to one of the $N$ modes, $f_1,\ldots,f_N$. A control trajectory is a \emph{schedule} composed of the timings and sequence of switches between the modes. A schedule of time length $t_f>0$ is given by $(\Sigma,\mathcal{T})$ where $\Sigma = [\sigma_1,\ldots,\sigma_M]$ is the mode sequence, $\mathcal{T} = [T_1,\ldots,T_{M-1}]$ is the vector of strictly monotonically increasing switching times, and $M$ is the number of modes in $\Sigma$. With $T_0 = 0$ and $T_M = t_f$, the dynamics at time $t$ where $T_{i-1}\leq t<T_i$, for some $i = 1,\ldots,M$ are
\[
\dot{x}(t) = f_{\sigma_i}(x(t)).
\]
We say the state and schedule pair $(x,(\Sigma,\mathcal{T}))$ is a dynamically \emph{feasible} trajectory if for all $t\in[0,t_f]$, there is an $i = 1,\ldots,M$ such that $T_{i-1}\leq t<T_i$ and
\begin{equation}
x(t) - x_0 - \sum_{j = 1}^{i-1}\int_{T_{j-1}}^{T_{j}}f_{\sigma_j}(x(\tau))d\tau -  \int_{T_{i-1}}^t f_{\sigma_i}(x(\tau))d\tau = 0.
\label{eq-S}
\end{equation}
Label the set of all such feasible pairs $(x,(\Sigma,\mathcal{T}))$ as $\mathcal{S}$.

The power network benchmark problem is the 118 Bus Test Case which has 118 busses, 54 of which are generators, and 26 switched capacitor banks which are all connected or disconnected in unison. The number of operating modes is $N = 2$.

\subsection{118 Bus Test Case}
The example power network has topology and line and bus parameters from the IEEE 118 Bus Test Case, a 1962 study of a segment of North America's midwest grid \citep{christie}. This network is composed of 118 buses, 186 lines, 54 generators and is shown in Fig.~\ref{fig-case118_Graph}.  The state of the system is the relative rotor angle and angular velocity of each generator excluding the reference generator, for a total of 106 states. 

In addition, we connect switched capacitor banks in series to 26 chosen transmission lines. The placement of the switched capacitors is a design decision discussed in Section~\ref{sec-design}. Each capacitor's capacitance is chosen so that when the switched capacitor is ``on'', its associated line's reactance doubles.  The location of each switched capacitor is shown in Fig.~\ref{fig-case118_Graph} and are chosen so that every generator is connected to at least one other adjacent generator through a line with a switched capacitor so that the $E_i$, $G_{ii}$, $Y_{ij}$, and $\psi_{ij}$ parameters can be switched for each generator's power output $P_{e,i}$, Eq.~\ref{eq-powerout}. For this study, all 26 switches are switched in unison so that all are ``on'' or ``off'' together. As such, the network has two modes of operation, $\dot{x}(t) = f_1(x(t))$ and $\dot{x}(t) = f_2(x(t))$.

The control is the scheduling of the fully coupled switching of all capacitor switches. We wish to schedule the capacitor switches such that network performance is improved. The improvement is measured by a cost function that assesses a greater performance to system trajectories for which the rotor's remain nearer steady state and rotor rotation frequency nearer operating frequency. The disturbance perturbation is a vector of random angles taken from a uniform distribution between $[-0.3,0.3]$ radians.

\section{Benchmark Problems}
\label{sec-problems}
The switched capacitor power network regulation problem is to reject a disturbance by driving the system toward synchronicity through connecting and disconnecting capacitors placed throughout the network. We introduce three problem types---\emph{mode scheduling}, \emph{sliding window scheduling}, and \emph{modal design}---and provide results for the benchmark problem

Optimal mode scheduling computes the timing of modal transitions and the sequence of transitions,  $(\Sigma, \mathcal{T})$, that optimizes a specified performance. It is an offline approach which computes the optimal state trajectory over a finite time interval. Sliding window scheduling is an online approach similar to model predictive control. It repeatedly computes the optimal (or near optimal) mode schedule for intervals of time that advance forward with real-time.  Modal design is the problem of designing the physical system to produce desirable modes of operation. For the benchmark problem the switched capacitor transmission line placement is a design choice.  Distinct placements correspond to distinct modes of operation which directly affects the network's effectiveness to reject a disturbance.

\subsection{Mode Scheduling \label{sec-scheduling}}
Mode scheduling has been thoroughly studied, but to the best of our knowledge, a preeminent strategy has not emerged.  The problem is to compute a feasible schedule $(\Sigma, \mathcal{T})\in\mathcal{S}$ that infimizes a cost function $J$. We only consider cost functions with the form 
\begin{equation}
J(x, (\Sigma, \mathcal{T})) = \int_0^{t_f}\ell(x(\tau))d\tau,
\label{eq-J}
\end{equation}
with running cost $\ell(x(\tau))\in\setR$ which does not directly rely on $(\Sigma, \mathcal{T})$. The mode scheduling problem is
\[
\inf_{(x, (\Sigma, \mathcal{T}))\in\mathcal{S}} J(x,(\Sigma, \mathcal{T}))
\]
The minimization is constrained to the set of dynamically feasible trajectories $\mathcal{S}$ (see Eq.~\ref{eq-S}). 

We apply the projection-based mode scheduling method we presented in \citep{caldwell_murphey_nahs} to solve for the optimal mode schedule (see Algorithm 1 with step size computation Algorithm 2 in \cite{caldwell_murphey_nahs}). The method is analogous to numerical optimization like steepest descent. It iteratively computes a descent direction, takes a step in the descending direction that satisfies a sufficient descent condition, and updates. The descent direction is the negative \emph{mode insertion gradient} which is the sensitivity of the cost $J$ to a switch in the schedule for infinitesimal duration \citep{axelsson_wardi_et_al,egerstedt_wardi_axelsson,caldwell_murphey_nahs,gonzalez_vasudevan_etal,wardi_egerstedt_hale,wardi_egerstedt_acc12}. Specifically, the mode $\sigma\in\{1,\ldots,N\}$ insertion gradient at time $t\in[0,t_f]$ is
\begin{equation}
\begin{array}{c}
d_\sigma(t;\Sigma,\mathcal{T}) := \rho(t)^T(f_\sigma(x(t)) - f_{\sigma_i}(x(t)))\textrm{, }\rho(t_f) = 0\\\textrm{where }T_{i-1}\leq t < T_i
\end{array}
\label{eq-mig}
\end{equation}
where $\Sigma = [\sigma_1,\ldots,\sigma_M]$, $\mathcal{T} = [T_1,\ldots, T_{M-1}]$ and the adjoint $\rho(t)$ is the solution to 
\[
\begin{array}{c}
\dot{\rho}(t) = -Df_{\sigma_i}(x(t))^T\rho(t)-D\ell(x(t))^T,\\\textrm{where }T_{i-1}<t<T_{i}.
\end{array}
\]
The step size is chosen through backtracking which computes a step size large enough to guarantee algorithm convergence. The update is a projection to $\mathcal{S}$. The solution to each step of the iteration is a dynamically feasible trajectory. Convergence is with respect to $\theta$, the minimal value of the mode insertion gradient $d$, Eq.~\ref{eq-mig}:
\begin{equation}
\theta(\Sigma,\mathcal{T}) := \min_{\sigma\in\{1,\ldots,N\},t\in[0,t_f]} d_{\sigma}(t;\Sigma,\mathcal{T}).
\label{eq-theta}
\end{equation}
The value $|\theta|$ is analogous to the norm of the gradient as the optimality condition in smooth optimization \citep{caldwell_murphey_nahs}.

We apply projection-based mode scheduling to schedule the capacitor switches in order to respond to the disturbance.  The cost is given by $\ell(x(t),u(t)) = 1/2 (\delta(t)-\bar{\delta}(t))^T(\delta(t)-\bar{\delta}(t)) + 1/40 (\dot{\delta}(t)-2\pi f_s)^T(\dot{\delta}(t)-2\pi f_s)$ where $\bar{\delta}(t)$ is the mean rotor phase at time $t$ and $f_s$ is the generator frequency.  The backtracking parameters are set to $\alpha = 0.4$ and $\beta = 0.1$. 

\begin{figure*}
\centering
\def\svgwidth{0.95\textwidth}
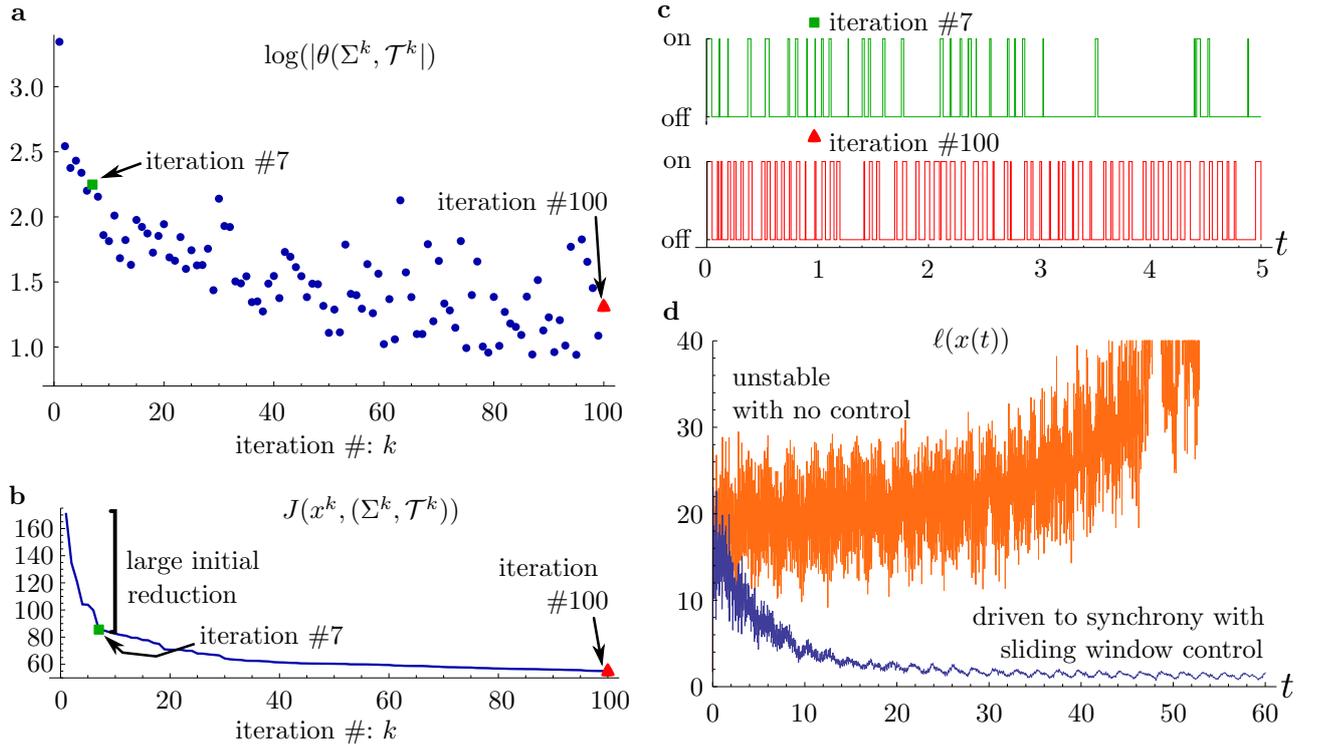
\caption{\textbf{a}, Convergence of optimality condition value $\theta$ toward zero as a function of iteration.  \textbf{b}, Cost $J$ as a function of iteration.  A large reduction occurs in the first 7 iterations. \textbf{c}, Comparison of the control signal for iterations 7 and 100. \textbf{d}, Comparison of running cost for no control (orange) and sliding window control (blue).}
\label{fig-power_net_results}
\end{figure*}

The results of mode scheduling the initial 5 seconds following a disturbance for 100 iterations of the algorithm are shown in Fig.~\ref{fig-power_net_results}. We find that the rotor phases do not diverge with the computed schedule. The cost reduces from $J = 170.68$ to $J = 54.78$ (see Fig.~\ref{fig-power_net_results}b), and the optimality function increases from $\theta = -2213.71$ to $\theta = -20.32$ (see Fig.~\ref{fig-power_net_results}a).  The total number of modes in the $7^\textrm{th}$ iteration's schedule is $M^7 = 66$, while the final switching schedule has $M^{100} = 120$. The schedules at the $7^\textrm{th}$ and $100^\textrm{th}$ iteration are in Fig.~\ref{fig-power_net_results}c.  

The shortest time between switches for any iteration is $10^{-8}$s. A direct mode scheduling method \citep{alamir_attia,lincoln_bernhardsson, lincoln_rantzer} would need to solve a $10^{10}-10^{11}$ dimensional optimization to capture equivalent fast switching over the 5s time horizon and 106 states.

For the initial iterations in which $(x^k,u^k)$ are far from an infima, both the optimality function (see Fig.~\ref{fig-power_net_results}a) and the cost (see Fig.~\ref{fig-power_net_results}b) reduce significantly, which is a phenomenon that commonly occurs with first-order smooth numerical optimization implementations like steepest descent.  Since most of the gained performance occurs in the first few iterations, it is reasonable to expect that a sliding window real-time approach is viable.  Such an approach computes only the first or first few mode scheduling iterations for each window.  

\subsection{Sliding Window Scheduling \label{sec-sliding}}
For online execution, the optimal schedule can be computed for each time window of a sliding window implementation. The model predictive control type approach applies the projection-based mode scheduling algorithm for intervals of time with fixed length that advance forward with time. In a real-time implementation, the advances in time occur with the update of system sensors. Specifically, the approach computes a schedule for a time window of duration $T = 5$ seconds but applies it for only $dt = 0.1$ seconds before incrementing the window $dt$ seconds and repeating for the new $5$ second time window.  The current window's initial state inherits the previous window's state at time $t_{i-1}+dt$. Instead of computing the full optimal schedule for each window, we compute a single projection-based mode scheduling iteration for the time interval $t\in[t_i,t_i+T]$. The goal is for a real-time active control rejection of the disturbance.

Fig.~\ref{fig-power_net_results}d compares the running cost $\ell(\cdot)$ for the sliding window result against the no control result.  Without control, the system destabilizes, while sliding window single-bit control drives the system toward synchrony. A core i7-3770K processor computes each window's schedule in an average of 1.94 seconds. While the current implementation is 20 times slower than real-time, it indicates that an improved implementation on a more advanced computing machine could execute the sliding window approach real-time even for the high-dimensional example.

\subsection{Modal Design \label{sec-design}}
Modal design is the problem of designing the physical system so that its evolution is dictated by the most desirable set of modes of operation possible. For the benchmark problem, the space of candidate designs is the switched capacitors' line placements. The designs are ranked by their capability to respond to a disturbance, where a disturbance is modeled as before as an initial perturbation from steady state. We quantify a design's capability through a performance metric's sensitivity to mode switches---recall the control authority is solely through mode switches.  In other words, the preferred design is the one for which the control authority can impose the greatest performance. 

\begin{figure}
\centering
\def\svgwidth{.44\textwidth}
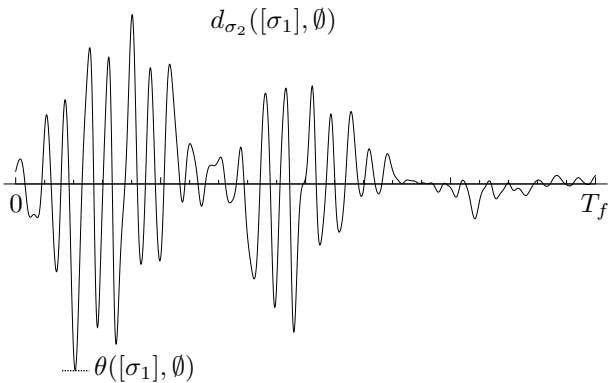
\caption{Example mode insertion gradient of mode $\sigma_2$ into schedule $([\sigma_1],\emptyset)$ with sensitivity value $\theta$ marked. }
\label{fig-bus39_MIG_sensitivity}
\end{figure}

Each distinct line placement design corresponds to a distinct set of modes of operation. Given a cost $J$, Eq.~\ref{eq-J}, a design's capability is a function of the metric's sensitivity to switches between the designed modes.  For schedule $(\Sigma,\mathcal{T})$, this sensitivity value is $\theta(\Sigma,\mathcal{T})$, Eq.~\ref{eq-theta}.  An example mode insertion gradient of mode $\sigma_2$ is in Fig.~\ref{fig-bus39_MIG_sensitivity} for the no switch mode schedule $\Sigma = [\sigma_1]$, $\mathcal{T} = \emptyset$, with $\theta$ marked.  

The value of the mode insertion gradient of a mode $\sigma_i$ at a time $t$ indicates the increase or decrease to the cost if $\sigma_i$ is inserted into the schedule for a sufficiently short duration at that time. When $d_{\sigma_i}(t;\Sigma,\mathcal{T})<0$, a sufficiently short duration insertion of $\sigma_i$ results in a reduced cost \citep{caldwell_murphey_nahs}. The optimal mode and timing arguments that compute $\theta$ in Eq.~\ref{eq-theta} correspond to the mode and timing for which an insertion for sufficiently short duration reduces the cost the most. 

The sensitivity value $\theta$ only provides an indication to the behavior of the cost for local schedule changes. Since a cost function's local behavior is usually a poor indicator of its optimal value, we find a correlation between the $\theta$ value for a large number of designs and the optimal cost for the benchmark problem. We selected 200 distinct sets of capacitor line placements with the requirement that for a placement set, every generator is connected to at least one other adjacent generator through a line with a switched capacitor. The process to generate a placement set was random. A random line for which at least one neighboring generator did not yet meet the requirement was repeatedly chosen until the requirement was fulfilled. Each capacitor's capacitance was chosen so that when the switched capacitors are ``on'', its associated line's reactance doubles. One could choose to have the capacitance as a design variable as well. As before, we consider two modes of operation $f_1$ and $f_2$, where $\dot{x} = f_1(x(t))$ when all capacitors are disconnected, or ``off'', and $\dot{x} = f_2(x(t))$ when all capacitors are connected, or ``on''.  

\begin{figure}
\centering
\def\svgwidth{.44\textwidth}
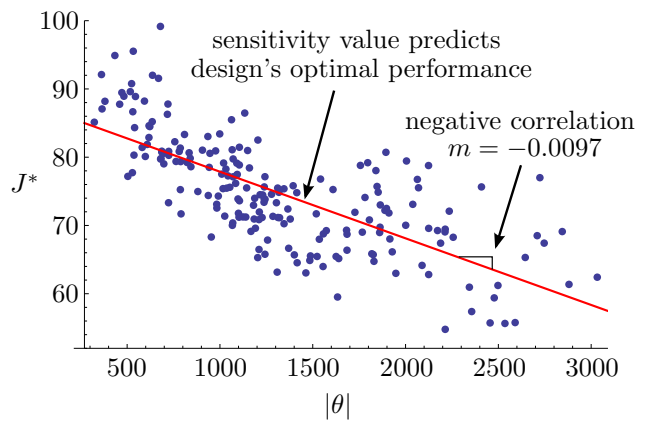
\caption{Correlation between the sensitivity value $\theta$ and the optimal cost $J^\star$ for the no switch schedule $([\sigma_1],\emptyset)$. }
\label{fig-optcost_vs_theta0}
\end{figure}

The correlation between the optimal cost $J^\star$, computed through projection-based mode scheduling (see Section~\ref{sec-scheduling}), and the sensitivity value $\theta$ for the no switch schedule $( [\sigma_1],\emptyset)$ is shown in Fig.~\ref{fig-optcost_vs_theta0}. We find that a design with a more negative $\theta$ is more likely to result in a lesser optimal cost and therefore be more capable at rejecting the disturbance.

The space of candidate designs is not limited to switched capacitor line placements. The candidate designs can include any set of controllable switching component with differing parameters as long as distinct designs correspond to distinct modes of operation.

\section{Discussion and Open Questions}
\label{sec-discussion}
We proposed a benchmark problem for switched-mode optimal control. The benchmark is the IEEE 118 Bus Test Case with dynamics given by the swing equations and control authority through connecting and disconnecting switched capacitor banks. The benchmark has 106 states for the 54 generating nodes (recall one generator acts as a reference) and 2 distinct operating modes for the coupled switching of all switched capacitors. 

Scaling in both state and the number of operating modes are important concerns. For example, can mode scheduling scale to networks with a greater number of generating nodes like North America's Western and Eastern Interconnections? The sliding window scheduling executed in about 20 times slower than real-time on a single processor. Could an improved implementation execute real-time for a larger network?

How does switched-mode optimal control scale to more operating modes? For the IEEE 118 Bus Test Case, only two operating modes were needed to respond to a disturbance as shown in Fig.~\ref{fig-power_net_results}d. At any time, only a single bit of control---i.e. ``on'' or ``off''---is needed, which lowers computational and communication complexity. Do other power networks require more operating modes? For the projection-based mode scheduling, the number of mode insertion gradients is $N$ but the values of each mode insertion gradient must be compared with each other, which at worst case grows $N^2$. 

The swing equations are a good representation, but are a simplification of real power network evolution, which brings into question the reliability of switched-mode solutions. Can sliding window scheduling or other feedback strategies reject model disturbances? How practical is sliding window scheduling to more complex, but accurate models, like that in \cite{natarajan_weiss}.

Another concern is that we assume that the switched capacitors can be switched arbitrarily fast, but there may be a significantly long transition period. Additionally, there may be preferred switching frequencies which would need to handle. 

\textbf{Open questions include:}
\begin{enumerate}
\item Can switched-mode optimal control scale in state? What about for real-time implementations?
\item How many operating modes, $N$, are needed to respond to a disturbance? Should every transmission line have independent switching control authority?
\item How reliable are switched-mode solutions? 
\item How does switched-mode optimal control scale to more complex models?
\item Can switching mode control efficiently implement hysteresis and timing constraints?
\end{enumerate}

Each question is directly relevant to the proposed benchmark problem making it an ideal test case for switched-mode optimal control assessment and development.

\section{Conclusion}
We proposed a benchmark problem for assessing and developing switched-mode optimal control approaches using the IEEE 118 Bus Test Case and the swing equations. Coupled with controllable switching components, the benchmark is a high-dimensional, nonlinear, switched-mode system. A baseline for relevant approaches like mode scheduling, sliding window scheduling, and modal design is given. Open questions include scalability in state and the number of modes of operation, as well as reliability, hysteresis, and timing constraints.

\balance

\bibliography{worksbib}             
\end{document}